\documentclass[10pt,letterpaper]{article}
\usepackage[top=0.85in,left=2.75in,footskip=0.75in,marginparwidth=2in]{geometry}

\usepackage[utf8]{inputenc}

\usepackage{cite}

\usepackage{nameref,hyperref}
\usepackage{amsmath}
\newtheorem{definition}{Definition}
\newtheorem{theorem}{Theorem}[section]

\newtheorem{lemma}[theorem]{Lemma}

\usepackage{amsfonts}
\usepackage{tikz}
\usetikzlibrary{shapes,arrows}

\usepackage[right]{lineno}

\usepackage{microtype}
\DisableLigatures[f]{encoding = *, family = * }

\raggedright
\usepackage{changepage}

\usepackage[aboveskip=1pt,labelfont=bf,labelsep=period,singlelinecheck=off]{caption}

\makeatletter
\renewcommand{\@biblabel}[1]{\quad#1.}
\makeatother

\usepackage{lastpage,fancyhdr,graphicx}
\usepackage{epstopdf}
\fancyheadoffset[L]{2.25in}
\fancyfootoffset[L]{2.25in}

\usepackage{color}

\definecolor{Gray}{gray}{.25}

\usepackage{graphicx}

\usepackage{wrapfig}
\usepackage[pscoord]{eso-pic}
\usepackage[fulladjust]{marginnote}
\reversemarginpar

\begin{document}
\vspace*{0.35in}

\begin{flushleft}
{\Large
\textbf\newline{A computational approximation for the solution of retarded functional differential equations and their applications to science and engineering}
}
\newline
\\
Burcu G\"{u}rb\"{u}z\textsuperscript{1}
\\
\bigskip
\bf{1} Institute of Mathematics, Johannes Gutenberg-University Mainz, Germany
\\
\bigskip
* burcu.gurbuz@uni-mainz.de

\end{flushleft}

\begin{abstract}
Delay differential equations are of great importance in science, engineering, medicine and biological models. These type of models include time delay phenomena which is helpful for characterising the real-world applications in machine learning, mechanics, economics, electrodynamics and so on. Besides, special classes of functional differential equations have been investigated in many researches. In this study, a numerical investigation of retarded type of these models together with initial conditions are introduced. The technique is based on a polynomial approach along with collocation points which maintains an approximated solutions to the problem. Besides, an error analysis of the approximate solutions is given. Accuracy of the method is shown by the results. Consequently, illustrative examples are considered and detailed analysis of the problem is acquired. Consequently, the future outlook is discussed in conclusion.
\end{abstract}
\bigskip
\section{Introduction}
\label{sec:1}
\bigskip
In general aspect, functional differential equations (FDEs) describe a wide range of real life phenomena which are extremely valuable in scientific research fields. Specifically, delay differential equations (DDEs) and their applications perform widely in many different branches of medicine, engineering and life sciences. These fields have operations based on  time-dependent processes which can be easily described by DDEs \cite{Bellen}, \cite{Halanay}.
\\~\\
Economic models are easily described by using  FDEs. Some arguments in the continuous systems are helpful to establish control problems and give us better understanding for the market models by using the time-delay contribution in equations. Multiple decision tasks, capable formulations are defined by FDEs with the idea of overlapping in time. Besides, competitive market scenarios are established with regard to the DDE models. In the theory of modern economic growth,  continuous-time modeling by FDEs derives an outstanding task. For instance, the Kaldor-Kalecki model, the Cournot oligopoly model and the Goodwin model are also well-known models in economy which are described by the differential equations with time-delay \cite{Keller}-\cite{Savku}.
\\~\\
In computer science, FDEs play an essential role. Recently, machine learning studies are very attractive in many applications. Time-delay differential equations establish a base for machine learning applications. Specifically, reservoir computing has been introduced recently as a machine learning pattern \cite{Grig1}. Its implementations with hardware accomplish processing of empirical data \cite{Grig2}. Besides, in machining processes chatter identification is described by FDEs \cite{Khasawneh}. In this area, the importance of FDEs has been described in current studies.
\\~\\
In electrodynamics, a device which is called as pantograph exists in an electric locomotive. This device collects the electric current and power is provided to the locomotive by a system from above cables. So, this system of electricity collection to trains or trams provides a contact with regard to the pantograph design \cite{Bhalekar}-\cite{Dehghan}. This system is a physical model in electrodynamics and it is related to the time dependence. The modeling of this structure is formulated by DDEs. In Fig. 1, the pantograph above a tram is shown.  This basic design of the pantograph is a physical motion and accomplishes DDEs with regard to the the vertical displacements of the pantograph at frame and head \cite{Abbott}-\cite{Gilbert}
\bigskip
\begin{figure}[htp]\label{f1}
\begin{center}
  \includegraphics[width=2in]{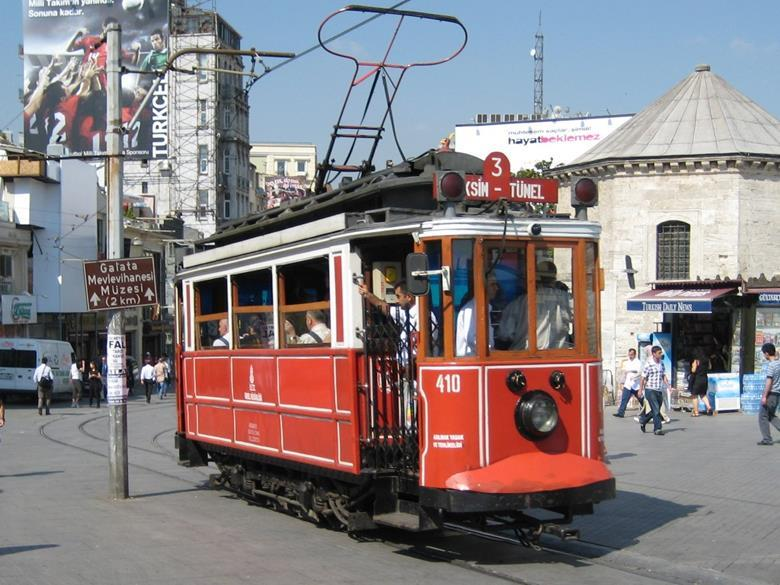}\\
  \caption{The pantograph of a tram (\.{I}stiklal Street, \.{I}stanbul) \cite{tram}.}\label{tram}
  \end{center}
\end{figure}
\bigskip
In the model, we can also consider some constants such as the damping coefficients for describing the motions from the head to the frame and other way around respectively with $\mu_{E}$ and $\mu_{D}$. A simple representation of a pantograph and its trolley wire system is shown in Fig. 2 which is formulated as
\bigskip
\\
\begin{eqnarray}\label{1}
 m_{E}y_{E}''(t)+\mu_{E}(y_{E}'(t)&-&y_{D}'(t))+k_{E}(y_{E}(t)-y_{D}(t))+\lambda_{y}(t)+m_{E}g=0, \\
 m_{D}y_{D}''(t)+\mu_{D}y_{D}'(t)&=&F_{u}+\mu_{D}(-y_{D}'(t)+y_{E}'(t))+k_{E}(-y_{D}(t)+y_{E}'(t))+m_{D}g.
\end{eqnarray}
\bigskip
\\~\\
Here, $k_{E}$ and $k_{D}$ values are the constants which are presenting bounce between frame and head of the pantograph, $m_{E}$ and $m_{D}$  symbolise the masses of head and frame, respectively, $y_{E}(t)$ and $y_{D}(t)$ are related to the displacement of head and frame, respectively, $\lambda_y$ is the vertical component force, $F_{u}$ is the upward force and $g$ is the gravity of earth \cite{Benet}.
\bigskip
\begin{figure}[htp]\label{f2}
\begin{center}
  \includegraphics[width=4in]{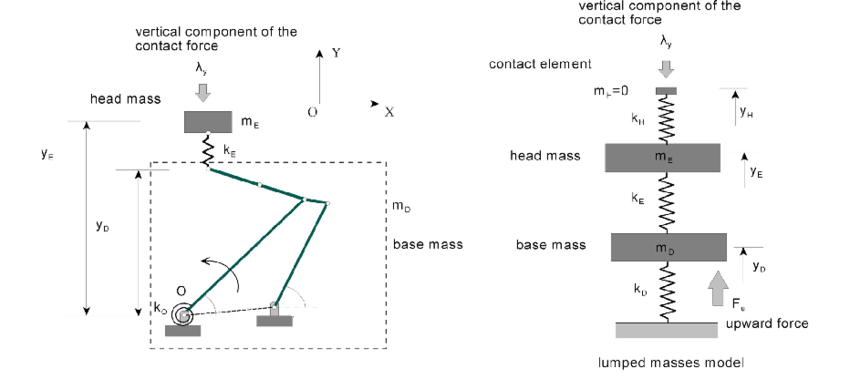}\\
  \caption{A representation of a pantograph and its trolley wire system \cite{Benet}.}\label{pantograph}
  \end{center}
\end{figure}
\bigskip
\\~\\
In biology, models which are described by the functional - differential equations have been established for determining a dynamics of some biological systems, and thermal science and mechanics systems. This wide usage of the FDEs supports mathematical biology studies and researchers in the field. For instance, recent studies on these type of models introduce a tool for explaining malaria, Dengue Fever epidemics and Covid-19 virus spreading models which include very frequently time delays \cite{Anderson}-\cite{Dell}. Furthermore, Mackey and Glass \cite{Mackey} proposed two possible models to describe the change of density of Hematopoietic cells in the blood that is circulating in the human body \cite{Bashier}.
\bigskip
\\~\\
Particularly, in population dynamics DDEs are often used. Continuous-time models including delay term describe mainly structure of a population. In this research field, the rate of producing egg in a population, susceptibility to parasitism or a bacterial population and their relations with related to the time dependence and delay term have been used in many studies \cite{Nisbet}-\cite{Diekmann}. Also, population dynamics is described by DDEs with the examples of investigation on baleen whale populations, delay model for introducing the dynamics of human immunodeficiency virus (HIV) infection etc. \cite{Nel}-\cite{Clark}.
\bigskip
\\~\\
Diversely, finding the solutions of these types of models analitically can be difficult. Due to this reason, numerical methods such as Adomian decompostion, finite element methods, homotopic perturbation techniques, Runge-Kutta type methods, direct two-point block methods, collocation methods, variational iteration methods have been applied to find approximate solutions of these type of problems \cite{Smith}-\cite{sb}.
\bigskip
\\~\\
In this study, a class of FDEs is described as following
\\
\begin{equation}\label{3}
u_{l}'(t)=-\gamma u_{l}(t)+ \beta u_{l}(t-\tau ) + g_{l}(t), \ 0 \leq t \leq b, \ l=1,2,3,
\end{equation}
\\
with initial condition
\\
\begin{equation}\label{4}
u_{l}(0)=\phi_{l},\ \ t \in [0,b].
\end{equation}
\\
Herein $\gamma, \beta >0$ are given as constants, $\tau>0$ is described for a delay term by using the name \textit{"time delay"}. Assuming that the approximate solutions of our problem, which has been given at (\ref{3})-(\ref{4}), defined in the form of Laguerre series as
\\
\begin{equation}\label{5}
\begin{aligned}
u_{l}(t)&\cong u_{l,N}(t)=\sum_{n=0}^{N}a_{l,n}L_{n}(t); \  \ 0 \leq t \leq b < \infty.
\end{aligned}
\end{equation}
\\
Here Laguerre polynomials are $L_{n}(t)$;
\\
\begin{equation}\label{6}
L_{n}(t)=\sum_{k=0}^{n}\frac{(-1)^k}{k!}\binom{n}{k}t^{k}, \ \ \  n \in N, \ \ \  0\leq t < \infty,
\end{equation}
\\
and $a_{n}$, $n=0,1,...,N$ are expressed the unknown coefficients of Laguerre polynomials. Moreover, $N>0$ is an integer which is selected properly, i.e. $N\geq 2$ \cite{Aizenshtadt}.
\\~\\
Applications of such dynamics in numerical studies are of valuable contribution to the field which enable us to understand the dynamical structure comprehensively. Here, an algorithmic approximation in the interdisciplinary concept is introduced which leads to widen related studies. The aim of the study is to give an alternative solution scheme and different perspective than the previous works in the field. In this way, it is aimed by this study that addresses the development of many fields through an important dynamic.
\\~\\
The organisation of the manuscript is given as follows. In Section 2, a mathematical outlook and preliminaries are given. In Section 3, the numerical method is defined with its details. Accuracy of the numerical technique is described in Section 4. In Section 5, numerical simulations are demonstrated by figures and tables. A brief conclusion and future studies are presented in Section 6.
\bigskip
\section{Preliminaries}
\label{sec:2}
\bigskip
Here some preliminaries about Laguerre polynomials and their properties on the delay differential equations has been introduced for the further applications. Besides, the recurrence relations of series and definitions with regard to the model are introduced to support the technical results.
\bigskip
\\~\\
Interdisiplinary researches of the Laguerre polynomials are well-known in many fields. The Laguerre polynomials are important in the applications of quantum theory in the hydrogen atom concept, chemical and mathematical physics \cite{Borwein}. They play an important role at the numerical methods in applied mathematics field such as Laguerre collocation method \cite{bg}, \cite{bgg}, \cite{bggg}, \cite{bgggg}, Laguerre spectral method \cite{Ben}-\cite{Baleanu}, Laguerre pseudospectral method \cite{Alici}, and Laguerre wavelets collocation method \cite{Zhou}. Besides, they support the improvement of the numerical solutions by combining with other numerical techniques, i.e. Laguerre-Galerkin methods \cite{Guo}, Laguerre Tau methods and so on \cite{Siyyam}.
\bigskip
\\~\\
Moreover, Laguerre polynomials are of a close relation with Hermite polynomials. Most of their properties and applications at the research area of orthogonal polynomials may help to find out characteristic of the Hermite polynomials. The Hermite polynomials $H_{n}(t)$ are connected to Laguerre polynomials as following:
\begin{align}\label{7}
        H_{2n}(t)   &=(-1)^{n}2^{2n}n!L_{n}^{(-1/2)}(t^2),\\
        H_{2n+1}(t) &=(-1)^{n}2^{2n+1}n!L_{n}^{(1/2)}(t^2).
\end{align}
In Fig. 3, we can see the relations between the first few Laguerre and Hermite polynomials with $n=0,1,2,3$ and $0\leq t\leq 2$ where Laguerre polynomials are presented by red lines and Hermite polynomials are shown by blue lines.
\bigskip
\begin{figure}[htp]\label{f3}
\begin{center}
  \includegraphics[width=2.5in]{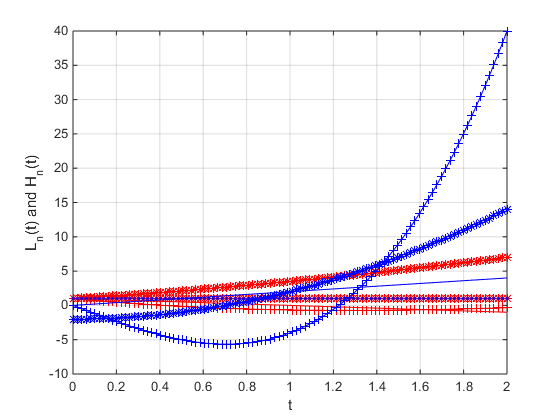}\\
  \caption{$L_n(t)$ and $H_n(t)$ values for $n=0,1,2,3$ and $t\in[0,2]$.}\label{LH}
  \end{center}
\end{figure}
\bigskip
\begin{definition} The first few Laguerre and Hermite polynomials are defined, respectively, as following:
\begin{align}\label{9}
&L_{0}(t)=1,                             &\ H_{0}(t)=1,\\
&L_{1}(t)=-t+1,                          &\ H_{1}(t)=2t,\\
&L_{2}(t)=\frac{1}{2}(t^2-4t+2),         &\ H_{2}(t)=4t^2-2,\\
&L_{3}(t)=\frac{1}{6}(-t^3+9t^2-18t+6).  &\ H_{3}(t)=8t^3-12t.
\end{align}
\end{definition}
Hermite polynomials have been applied on physical models, more specifically, at the solution of simple harmonic oscillator of quantum mechanics and Laguerre polynomials are seen in wave functions of hydrogen atom \cite{Arfken}. Usage of Laguerre polynomials are also well defined for the series solution. In this study, Laguerre polynomials give us an alternative series solution for the delay differential equations which is supported by following lemma.
\begin{lemma} \label{L: Laguerre} Suppose that $u(t)\in C^{(q)}[0,\infty[$ and $u^{(q+1)}$ is defined on the interval $[0,\infty[$ which is a piecewise continuous function. Then for
\begin{equation}\label{13}
 \mathcal{L}u(t)=\frac{d^q}{dt^q}u(t),
\end{equation}
we have converges uniformly $\sum_{i=0}^{\infty}u_{i}^{q}L_{i}(t)$ on $[0,\infty[$ to $\mathcal{L}u(t)$ for any positive integer $q$ \cite{Siyyam}.
\end{lemma}
In this study, a special class of functional differential equations has been considered. Functional differential equations are of many different types including delay differential equations or more specifically retarded, neutral and advanced functional differential equations \cite{Bellman}. Now, the definitions with regard to these equations are introduced in order to complete the problem settlement.
\begin{definition}
Here the following equation is considered:
\bigskip
\begin{equation}\label{14}
\zeta u'(t)+\gamma u(t)+\beta u(t-\tau)=g(t),
\end{equation}
which is given with $1$st order in its derivative and difference. In this place, an equation in the form above is called as \textit{retarded type} where $\zeta \neq 0$.
\end{definition}

\bigskip
\section{Method}
\label{sec:3}
\bigskip
In this section, a general scheme with regard to Laguerre collocation technique is introduced. We provide a novel collocation structure by using the Laguerre polynomial approach directly with related to some fundamental relations between the matrices. Here, the idea of the implementation of the technique is presented together with collocation points. The truncation of the series is shown and algebraic equations are obtained. We have the series approach and show the approximate solutions with the help of these connections.
\bigskip
\\~\\
Our motivation here is to show that we expand the well-conditioned collocation method with regard to Laguerre matrix approach. The applications of the technique some real-world models give us a chance to investigate the models based on numerical concept.
\bigskip
\subsection{Fundamental Relations}
\label{subsec:1}
\bigskip
In this subsection, an approximate solution with related to Laguerre polynomials form (\ref{6}) is presented. For our purpose, we consider a matrix form of Eq. (\ref{5}) as
\begin{equation}\label{15}
[ u_{l}(t)]=\mathbf{L}(t)\mathbf{A}_{l},
\end{equation}
in which
\begin{eqnarray}\label{16}
\begin{split}
\mathbf{L}(t)&=\left[\begin{array}{cccc}
    L_{0}(t) & L_{1}(t) & \cdots & L_{N}(t) \\
  \end{array}
  \right], \ \ \ \textrm{and}\\
  \mathbf{A}_{l}&=\left[
   \begin{array}{cccc}
     a_{i,0} & a_{i,1} & \cdots &a_{i,N}  \\
   \end{array}
 \right]^{T}, \ \ \  i=0,1,\dots,N.
       \end{split}
\end{eqnarray}
Then, we use the matrix relation
\begin{equation} \label{17}
 \mathbf{L}(t)=\mathbf{X}(t)\mathbf{H},
\end{equation}
where
\begin{eqnarray}\label{18}
\begin{split}
\mathbf{X}(t)&=\left[\begin{array}{cccc}
    1 & t & \cdots & t^{N} \\
  \end{array}
  \right],\\
 \mathbf{H}&=\left[
      \begin{array}{ccccc}
        \frac{(-1)^0}{0!}\left(
                     \begin{array}{c}
                       0 \\
                       0 \\
                     \end{array}
                   \right)
         &         \frac{(-1)^0}{0!} \left(
                     \begin{array}{c}
                       1 \\
                       0 \\
                     \end{array}
                   \right)
         &  \cdots
         &  \frac{(-1)^0}{0!}\left(
                     \begin{array}{c}
                       N \\
                       0 \\
                     \end{array}
                   \right) \\
0
        &\frac{(-1)^1}{1!} \left(
                     \begin{array}{c}
                       1 \\
                       1 \\
                     \end{array}
                   \right)
        &  \cdots
        & \frac{(-1)^1}{1!} \left(
                     \begin{array}{c}
                       N \\
                       1 \\
                     \end{array}
                   \right) \\
        \vdots &  \vdots & \ddots & \vdots \\
        0
        & 0
        & \cdots
        & \frac{(-1)^N}{N!}\left(
                     \begin{array}{c}
                       N \\
                       N \\
                     \end{array}
                   \right) \\
      \end{array}
    \right].
      \end{split}
\end{eqnarray}
Besides, the connections are defined between $\mathbf{X}(t)$ and $\mathbf{X}'(t)$. So, we write them as
\begin{equation}\label{19}
\mathbf{X}'(t)=\mathbf{X}(t)\mathbf{B},
\end{equation}
where
\begin{eqnarray}\label{20}
 \mathbf{B}=\left[
    \begin{array}{ccccc}
      0 & 1 & 0 & \cdots & 0 \\
      0 & 0 & 2 & \cdots & 0 \\
      \vdots & \vdots & \vdots & \ddots & \vdots \\
      0 & 0 & 0 & \cdots & N \\
      0 & 0 & 0 & \cdots & 0 \\
    \end{array}
  \right].
\end{eqnarray}
Then from the relations in (\ref{15}) and (\ref{19}), we get
\begin{equation}\label{21}
\mathbf{L}'(t)=\mathbf{X}(t)\mathbf{B}\mathbf{H},
\end{equation}
So that, from (\ref{17}), (\ref{19}) and (\ref{21}), we have
\begin{equation}\label{22}
\mathbf{X}(t)=\mathbf{L}(t)\mathbf{H}^{-1},
\end{equation}
\begin{equation}\label{23}
\mathbf{L}'(t)=\mathbf{L}(t)\mathbf{C},
\end{equation}
where
\begin{equation}\label{24}
\mathbf{C}=[c_{pq}], \ \
c_{pq}=\begin{cases}
      -1, & \text{if}\ p<q \\
      0, & \text{if}\ p\geq q.
    \end{cases}
\end{equation}
Then we have
\begin{equation}\label{25}
[u_{l}'(t)]=\mathbf{L}(t)\mathbf{C}\mathbf{A}_{l},
\end{equation}
Moreover, $t\rightarrow t-\tau$ is replaced into (\ref{15}) and get
\begin{equation}\label{26}
[u_{l}(t-\tau)]=\mathbf{L}(t-\tau)\mathbf{A}_{l}=\mathbf{X}(t)\mathbf{T}(t-\tau)\mathbf{B}\mathbf{H}\mathbf{A}_{l},
\end{equation}
where
\begin{eqnarray}\label{27}
 \mathbf{T}(t-\tau)=\left[
      \begin{array}{ccccc}
                \left(
                     \begin{array}{c}
                       0 \\
                       0 \\
                     \end{array}
                \right)
                     (-\tau)^0

         &         \left(
                     \begin{array}{c}
                       1 \\
                       0 \\
                     \end{array}
                   \right)(-\tau)^1
         &  \cdots
         &          \left(
                     \begin{array}{c}
                       N \\
                       0 \\
                     \end{array}
                   \right)(-\tau)^N \\
0
        &           \left(
                     \begin{array}{c}
                       1 \\
                       1 \\
                     \end{array}
                   \right)(-\tau)^0
        &  \cdots
        &           \left(
                     \begin{array}{c}
                       N \\
                       1 \\
                     \end{array}
                   \right) (-\tau)^{N-1} \\
        \vdots &  \vdots & \ddots & \vdots \\
            0
        &   0
        & \cdots
        &           \left(
                     \begin{array}{c}
                       N \\
                       N \\
                     \end{array}
                   \right)(-\tau)^0 \\
      \end{array}
    \right].
\end{eqnarray}
Then
\begin{eqnarray}\label{28}
   \mathbf{u}(t)&=& \overline{\mathbf{L}}(t)\mathbf{A}, \\
   \mathbf{u}'(t)&=& \overline{\mathbf{L}}(t)\overline{\mathbf{C}}\mathbf{A}, \\
   \mathbf{u}(t-\tau)&=&\overline{\mathbf{X}}(t)\overline{\mathbf{T}}(t-\tau)\overline{\mathbf{B}}\overline{\mathbf{H}}\mathbf{A},
\end{eqnarray}
are defined and in which
\begin{eqnarray}\label{31}
\mathbf{u}'(t)&=&\left[
              \begin{array}{c}
                \mathbf{u}'_{1}(t) \\
                \mathbf{u}'_{2}(t) \\
                \mathbf{u}'_{3}(t) \\
              \end{array}
      \right],\
\mathbf{u}(t-\tau)=\left[
              \begin{array}{c}
                \mathbf{u}_{1}(t-\tau) \\
                \mathbf{u}_{2}(t-\tau) \\
                \mathbf{u}_{3}(t-\tau) \\
              \end{array}
      \right],\
\overline{\mathbf{L}}(t)=\left[
              \begin{array}{ccc}
                \mathbf{L}(t)& 0  & 0  \\
                0 & \mathbf{L}(t)   & 0 \\
                0 & 0 &  \mathbf{L}(t) \\
              \end{array}
      \right], \\
\overline{\mathbf{X}}(t)&=&\left[
              \begin{array}{ccc}
                \mathbf{X}(t)& 0  & 0  \\
                0 & \mathbf{X}(t)   & 0 \\
                0 & 0 &  \mathbf{X}(t) \\
              \end{array}
      \right], \
\overline{\mathbf{B}}=\left[
              \begin{array}{ccc}
                \mathbf{B}& 0  & 0  \\
                0 & \mathbf{B}& 0 \\
                0 & 0 &  \mathbf{B} \\
              \end{array}
      \right], \
\overline{\mathbf{H}}=\left[
              \begin{array}{ccc}
                \mathbf{H}& 0  & 0  \\
                0 & \mathbf{H}& 0 \\
                0 & 0 &  \mathbf{H} \\
              \end{array}
      \right],\\
\overline{\mathbf{T}}(t-\tau)&=&\left[
              \begin{array}{ccc}
                 \mathbf{T}(t-\tau)& 0  & 0  \\
                0 &  \mathbf{T}(t-\tau)& 0 \\
                0 & 0 &   \mathbf{T}(t-\tau) \\
              \end{array}
      \right],\
\overline{\mathbf{C}}=\left[
              \begin{array}{ccc}
                \mathbf{C}& 0  & 0  \\
                0 & \mathbf{C}& 0 \\
                0 & 0 &  \mathbf{C} \\
              \end{array}
      \right],\
\mathbf{A}=\left[
              \begin{array}{c}
                \mathbf{A}_{1} \\
                \mathbf{A}_{2} \\
                \mathbf{A}_{3} \\
              \end{array}
      \right].
\end{eqnarray}
By using the equations (\ref{28})-(30), a matrix form of (\ref{3}) is described as
\begin{equation}\label{34}
\{\overline{\mathbf{L}}(t)\overline{\mathbf{C}}+\gamma \overline{\mathbf{L}}(t)-\beta \mathbf{u}(t-\tau)\}\mathbf{A}=\mathbf{g}(t).
\end{equation}
Here
\begin{equation}\label{35}
  \mathbf{g}(t)=\left[
              \begin{array}{c}
                g_{1}(t) \\
                g_{2}(t) \\
                g_{3}(t) \\
              \end{array}
      \right].
\end{equation}
\bigskip
\subsection{Method of Solution}
\label{subsec:2}
\bigskip
In this subsection, collocation method is applied by using the collocation points which give us a pointwise approximation with regard to our truncation number $N$ together with step size $h$. Then by replacing the collocation points
\begin{equation}\label{36}
 t_{i}=\frac{s}{N}i, \ \ \textrm{where}  i=0,1,\dots,N, \ \ \   \textrm{and} \ \ \  h=\frac{i}{N},
\end{equation}
into Eq. (\ref{34}), so that the fundamental matrix equation is obtained and written as follows:
\begin{equation}\label{37}
\{\overline{\mathbf{L}}(t_{i})\overline{\mathbf{C}}+\gamma \overline{\mathbf{L}}(t_{i})-\beta \mathbf{u}(t_{i}-\tau)\}\mathbf{A}=\mathbf{g}(t_{i})=\mathbf{G}.
\end{equation}
Briefly,
\begin{equation}\label{38}
\mathbf{W}(t_{i})\mathbf{A}=\mathbf{G}; \ \ [\mathbf{W};\mathbf{G}].
\end{equation}
Correspondingly, the initial condition is written in the matrix form as
\begin{equation}\label{39}
[u_{l}(0)]=\mathbf{L}(0)\mathbf{C}\mathbf{A}_{l}= [\phi_{l}],
\end{equation}
or
\begin{equation}\label{40}
\mathbf{U}\mathbf{A}_{l}= \mathbf{\phi}_{l}.
\end{equation}
Therefore, the last rows of the augmented matrix (\ref{38}) is reestablished by the row matrix (\ref{39}) in order to get an approximate solution of Eq.  (\ref{3}) under the condition (\ref{4}).  Thus, the augmented matrix is obtained as
\begin{equation}\label{41}
[\tilde{\mathbf{W}};\tilde{\mathbf{G}}].
\end{equation}
Laguerre coefficients can be computed by solving the augmented matrix system \cite{Cetin}-\cite{Gulsu}. Hence, the approximate solution is established as in the desired form of Laguerre series (\ref{6}) as
\begin{equation}\label{42}
\begin{aligned}
u_{l}(t)&\cong u_{l,N}(t)=\sum_{n=0}^{N}a_{l,n}L_{n}(t); \  \ 0 \leq t \leq b < \infty , \ \  l=1,2,3.
\end{aligned}
\end{equation}
\bigskip
\section{Accuracy}
\label{sec:4}
\bigskip
In this section, errors of approximation and solution strategy with regard to the its algorithm have been introduced. Error analysis has been proposed by using different norms in order to show the results more comprehensively. This lead us to see efficiency of our numerical technique. Besides, algorithm has been presented which is based on programming part of the construction of our method. Hereby the application of the numerical method by using a computer programming has been shared.
\bigskip
\\~\\
Accuracy is an important topic in approximation theory which has a reason to be introduced due to the errors occur in rounding and truncation of numerical solutions. When we use numerical algorithms with finite sensitivity, approximation errors, round or truncation error at the series solutions appears. Due to this reason, we should show the accuracy of the technique and how we deal with this approximation at the end of our application.
\bigskip
\subsection{Error Estimation}
\label{subsec:3}
\bigskip
In this subsection, a short presentation of the error estimation for the solutions with regard to Laguerre approach (\ref{6}) is given. Hereby accuracy of this technique is supported by this approximation. Let us first describe, $E_{N}(t)$, \textit{"error function" }for $t=t_{\alpha}, \alpha=0,1,...$
\begin{equation}\label{43}
\begin{aligned}
E_{l,N}(t_{\alpha})=|u_{l}(t_{\alpha})-u_{l}'(t_{\alpha})-\gamma u_{l}(t_{\alpha})+ \beta u_{l}(t_{\alpha}-\tau ) + g_{l}(t_{\alpha})|\cong0,
\end{aligned}
\end{equation}
where $k$ is appointed and we pay attention to the truncation limit $N$ which is an important argument for finding the approximate results. Accordingly, the difference $E_{N}(t_{\alpha})$ becomes smaller and, in here, this is denoted by $10^{-k}$. Diversely, some particular norms are given in order to investigate the distinctive error functions for measuring errors. So, these are described as follows.
\begin{itemize}
  \item For $L_2$ ; $E_{l,N}(t_{\alpha})=(\sum_{i=1}^{n}(e_{i})^2)^{1/2}$,
  \item For $L_{\infty}$ ; $E_{l,N}(t_{\alpha})=Max(e_{i})$ for $0\leq i\leq n$,
  \item For $RMS$ ; $E_{l,N}(t_{\alpha})=\sqrt{\frac{\sum_{i=1}^{n}(e_{i})^2}{n+1}}$.
\end{itemize}
Here $e_{l}=u(t_{l})-u_{N}(t_{l})$ also $u_{N}(t_{l})$ and $u(t_{l})$ are approximate and exact solutions of the problem, respectively.
\subsection{Algorithm}
\label{subsec:4}
\bigskip
In this section, \textit{Steps} of algorithm for the present method:
\begin{enumerate}
  \item[ ] \textbf{Data:} $\gamma$, $\beta$ and $\tau$ constants in Eq. (\ref{3}).
  \item[ ] \textbf{Result:} $u_{l,N}(t)$: approximate solutions.
  \item[\textit{S0.}] Truncation is chosen as $n\leq N$ for $n \in \mathbb{N}$,
  \item[\textit{S1.}] Construction of all the matrices,
  \item[\textit{S2.}] Replacement of the fundamental matrix equation,
  \item[\textit{S3.}] Apply the collocation points (colloc. pts.), $t_{i}=\frac{b}{N}i, \  i=0,1,\dots,N$ and $h=\frac{i}{N}$, to the fundamental matrix equation in \textit{S2}.
  \item[\textit{S4.}] Computation of the augmented matrix $[\mathbf{W}:\mathbf{U};\mathbf{G}]$ by Gauss elimination,
  \item[\textit{S5.}] Construction of the initial conditions (ICs) in matrix forms $[u_{l}(0)]$,
  \item[\textit{S6.}] Replacement of the initial condition in matrix forms in \textit{S5}. to the augmented matrix in \textit{S4}. Then we get $[\tilde{\mathbf{W}}:\tilde{\mathbf{U}};\tilde{\mathbf{G}}]$,
  \item[\textit{S7.}] Solution of the system in \textit{S6}. and replacement in the truncated Laguerre series form in Eq. (\ref{5}).
  \item[\textit{S8.}] Stop.
\end{enumerate}
We found an effective algorithmic approach for calculating an approximation procedure and to investigate the dynamics behind our model. This made a considerable impact on the dynamics and it has an excellent significance framework on approximation methodology.
\\~\\
Algorithms particularly connect the applicability with time standard. In order to increase the quality and to decrease the central processing unit (CPU) time, we may use an effective approach for coding. Due to this reason, straightforward techniques are of usable and effective to get the approximations. Then we can see in \ref{subsec:4} that our numerical method may be applied with regard to this idea. Besides, the following flowchart of the algorithm shows us a standardized approach to find the approximate solutions in our problem (\ref{3})-(\ref{4}) \cite{rum}.
\\
\tikzstyle{decision} = [diamond, draw, fill=blue!20,
    text width=4.5em, text badly centered, node distance=3cm, inner sep=0pt]
\tikzstyle{block} = [rectangle, draw, fill=blue!20,
    text width=5em, text centered, rounded corners, minimum height=4em]
\tikzstyle{line} = [draw, -latex']
\tikzstyle{cloud} = [draw, ellipse,fill=red!20, node distance=3cm,
    minimum height=2em]

\begin{tikzpicture} \label{fc}[node distance = 2cm, auto]
    \node [block] (init) {Initialize model};
    \node [cloud, left of=init, node distance=3cm] (Input) {Input: $\gamma$, $\beta$ and $\tau$};
    \node [cloud, right of=init, node distance=3cm] (Truncation) {Truncation: $n$};
    \node [block, below of=init, node distance=2cm] (identify) {Identify related matrices};
    \node [block, below of=identify, node distance=2cm] (evaluate) {Evaluate fundamental matrix equation};
    \node [cloud, right of=evaluate, node distance=3cm] (Collocation) {Colloc. pts. $t_{i}$};
    \node [block, below of=evaluate, node distance=2cm] (Augmented) {Compute augmented matrix};
    \node [cloud, right of=Augmented, node distance=3cm] (Condition) {ICs: $u_{l}(0)$};
    \node [block, below of=Augmented, node distance=2cm] (NewAugmented) {New augmented matrix};
    \node [block, left of=Augmented, node distance=3cm] (update) {Update truncation};
    \node [block, below of=NewAugmented, node distance=2cm] (Solve) {Solve the system};
    \node [decision, below of=Solve, node distance=2cm] (decide) {Is it the best candidate solution?};
    \node [block, below of=decide, node distance=3cm] (stop) {Stop};
    \path [line] (init) -- (identify);
    \path [line] (identify) -- (evaluate);
    \path [line] (evaluate) -- (Augmented);
    \path [line] (Augmented) -- (NewAugmented);
    \path [line] (NewAugmented) -- (Solve);
    \path [line] (Solve) -- (decide);
    \path [line] (decide) -| node [near start] {No} (update);
    \path [line] (update) |- (identify);
    \path [line] (decide) -- node {Yes}(stop);
    \path [line,dashed] (Input) -- (init);
    \path [line,dashed] (Truncation) -- (init);
    \path [line,dashed] (Truncation) -- (Collocation);
    \path [line,dashed] (Collocation) -- (evaluate);
    \path [line,dashed] (Truncation) |- (identify);
    \path [line,dashed] (Condition) -- (Augmented);
\end{tikzpicture}
\bigskip
\section{Numerical Illustrations}
\label{sec:5}
\bigskip
In here, some illustrative numerical examples are presented to show applicability of this given present method. These examples have been chosen from science and engineering applications to be analysed. Besides, the algorithm, which has been explained \ref{subsec:4}, is applied on the problems. The problem has been introduced with the general formula at (\ref{3})-(\ref{4}). Maple and Matlab computer programmes have been used for the calculations for finding the results and plotting the figures \cite{Maple}-\cite{MATLAB}.
\bigskip
\\~\\
\textbf{Example 1.}
\bigskip
Let us first consider a model from biology which is called as \textit{"Wazewska-Czyzewska and Lasota model"}. This biological model describes a remainder red blood cells, specifically, in animals \cite{Su}, \cite{Wazewska}:
\begin{equation}\label{44}
\frac{du(t)}{dt}=-\gamma u(t)+ \beta \exp (-\rho u(t-\tau )) , \  t > \tau,
\end{equation}
and its initial condition
\begin{equation}\label{45}
u(t)=\sin (t) \ \textrm{for} \ 0\leq t\leq \tau.
\end{equation}
Here, number of total red blood cells at time, denoted by $t$, is expressed by $u(t)$. Moreover, $\gamma$ which is presenting death rate of a unique red blood cell. Besides, $\beta$ and $\rho$ are denoted red blood cells generation per unit time and $\tau$ is time essential for generating a red blood cell. Specifically, $l=1$, $\gamma=0.4$, $\beta=\rho=1$ are chosen. So that, we have (28) as
\begin{equation}\label{46}
\frac{du(t)}{dt}=-0.4 u(t)+  \exp (-u(t-\tau )), \ \textrm{where} \  t > \tau .
\end{equation}
Then we have the matrix form as
\begin{equation}\label{47}
(\mathbf{L}(t)\mathbf{C}+0.4 \mathbf{L}(t)-\mathbf{X}(t)\mathbf{T}(t - \tau)\mathbf{B}\mathbf{H})\mathbf{A}=[g(t)].
\end{equation}
By using the procedure, we obtain approximate solutions. The approximate solutions for different $N$ values can be seen in Figure 1. Different error norm results for $N=3$ is given by Table 1.
\bigskip
\begin{figure}[htp]\label{f4}
\begin{center}
  \includegraphics[width=2in]{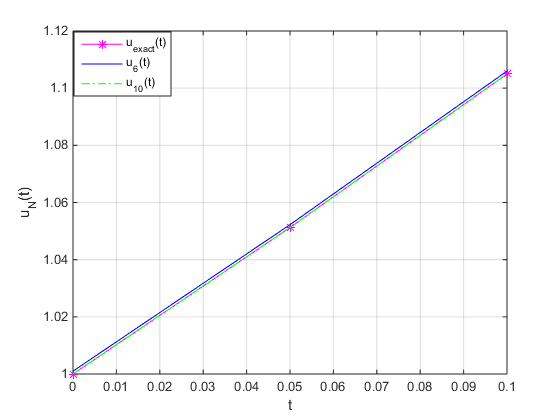}\\
  \caption{Exact solution comparison with some approximate solutions for $N=6$ and $10$.}\label{dcdss1}
  \end{center}
\end{figure}
\bigskip
\bigskip
\begin{table}[htp]
\centering
\caption{$L_\infty$, $L_2$ and $RMS$ errors for $N=3$.}
\label{tab:1}       
\begin{tabular}{llll}
\hline\noalign{\smallskip}
$t$ & $L_2$-Error & $L_\infty$-Error & $RMS$-Error \\
\noalign{\smallskip}\hline\noalign{\smallskip}
1 & 0.7560E-05 & 0.5247E-04 & 0.1000E-06\\
2 & 0.1164E-05 & 0.3791E-04 & 0.1502E-05\\
3 & 0.1550E-04 & 0.5467E-03 & 0.6855E-04\\
4 & 0.8259E-03 & 0.7795E-03 & 0.1752E-03\\
5 & 0.4643E-04 & 0.5467E-02 & 0.2916E-05\\
\noalign{\smallskip}\hline
\end{tabular}
\end{table}
\bigskip
\textbf{Example 2.}
\bigskip
Here, the model is considered as in the form \cite{Jumaa}, \cite{Shampine}:
\begin{eqnarray}\label{48}
\frac{du_{1}(t)}{dt}&=& u_{1}(t-2), \  t \geq 0, \\
\frac{du_{2}(t)}{dt}&=& u_{1}(t-2)+u_{2}(t-0.5), \  t \geq 0,
\end{eqnarray}
with initial conditions
\begin{equation}\label{49}
u_{1}(t)=1 \ \ \textrm{and} \ \ u_{2}(t)=1, \  t \leq 0.
\end{equation}
Herein, $l=2$ is chosen. For Eq. (39) $\gamma=0$, $\beta=1$, $\tau=2$ and for Eq. (40) $\gamma=0$, $\beta=1$ are assigned. Moreover, for $u_{1}(t)$, $\tau=2$ and for $u_{2}(t)$, $\tau=0.5$ are defined in Eq. (40). Again, we follow the similar process by implementing Laguerre collocation method (LCM) for $N=3$ and $N=4$ truncation values on our problem. Besides, this results obtained by LCM, some other numerical methods such as Runge-Kutta method (RKM) with fourth-order and Hermite collocation method (HCM) with the same $N$ truncation values. The results can be seen by the figures and the tables.
\bigskip
\begin{figure}[htp]\label{f5}
\begin{center}
  \includegraphics[width=2in]{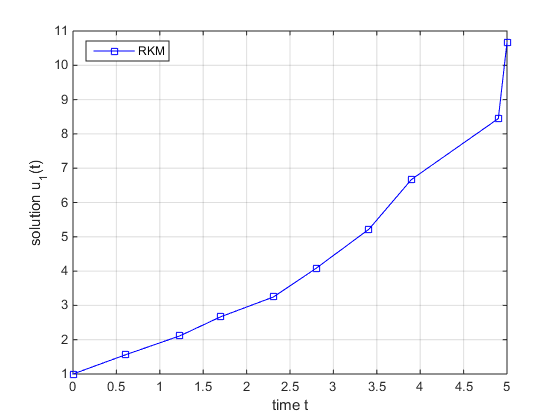}\\
  \caption{Runge-Kutta method (RKM) solution for $u_{1}(t)$ of Example 2.}\label{inbak1_21}
  \end{center}
\end{figure}
\bigskip
\begin{figure}[htp]\label{f6}
\begin{center}
  \includegraphics[width=2in]{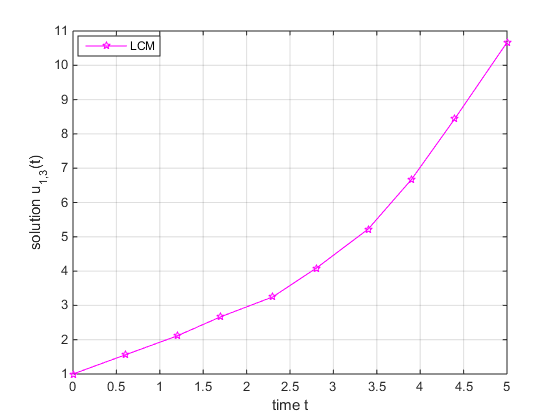}\\
  \caption{Laguerre collocation method (LCM) solution for $u_{1}(t)$, $N=3$ of Example 2.}\label{inbak1_22}
  \end{center}
\end{figure}
\bigskip
\bigskip
\begin{figure}[htp]\label{f7}
\begin{center}
  \includegraphics[width=2in]{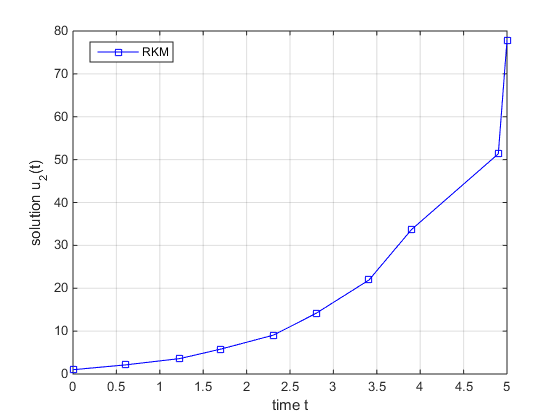}\\
  \caption{Runge-Kutta method (RKM) solution for $u_{2}(t)$ of Example 2.}\label{inbak1_23}
  \end{center}
\end{figure}
\bigskip
\bigskip
\begin{figure}[htp]\label{f8}
\begin{center}
  \includegraphics[width=2in]{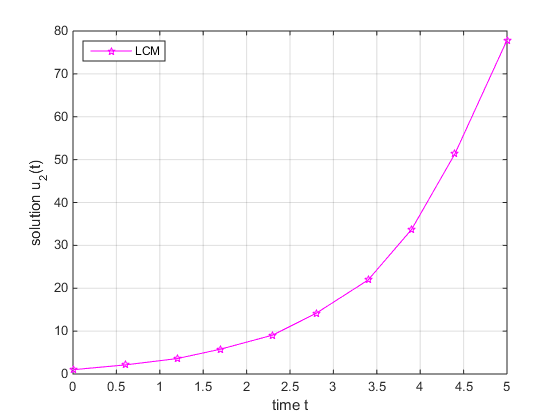}\\
  \caption{Laguerre collocation method (LCM) solution for $u_{2}(t)$, $N=3$ of Example 2.}\label{inbak1_24}
  \end{center}
\end{figure}
\bigskip
\bigskip
\begin{figure}[htp]\label{f9}
\begin{center}
  \includegraphics[width=2in]{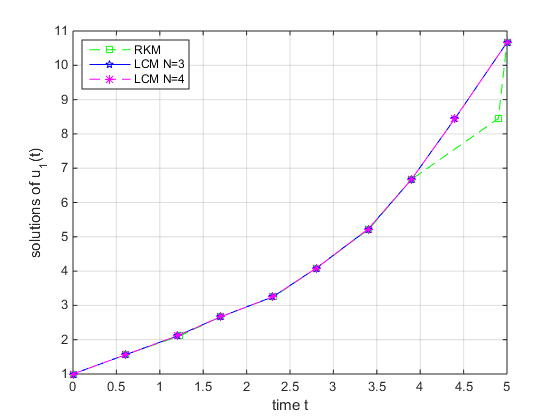}\\
  \caption{Comparison between RKM and LCM solutions for $N=3, 4$ and $u_{1}(t)$ of Example 2.}\label{inbak1_23}
  \end{center}
\end{figure}
\bigskip
\bigskip
\begin{figure}[htp]\label{f10}
\begin{center}
  \includegraphics[width=2in]{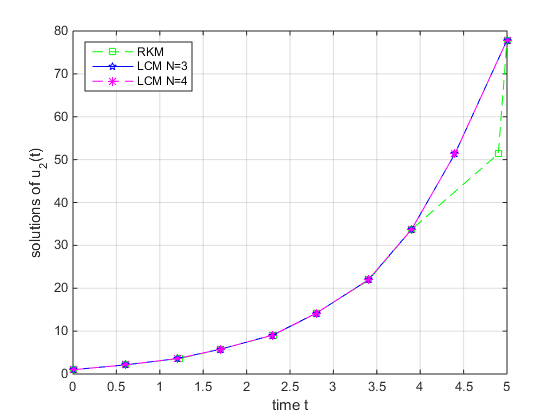}\\
  \caption{Comparison between RKM and LCM solutions for $N=3, 4$ and $u_{2}(t)$ of Example 2.}\label{inbak1_23}
  \end{center}
\end{figure}
\bigskip
\bigskip
\begin{table}[htp]
\centering
\caption{$L_2$ errors of LCM and HCM for $N=3$ and $N=4$ of Example 2.}
\label{tab:2}       
\begin{tabular}{lllll}
\hline\noalign{\smallskip}
$t$ & LCM, $N=3$ & LCM, $N=4$ & HCM, $N=3$ & HCM, $N=4$\\
\noalign{\smallskip}\hline\noalign{\smallskip}
0.0 & 0.6250E-02 & 0.1593E-03 & 0.7081E-02  & 0.2301E-03\\
1.2 & 0.4202E-03 & 0.2490E-04 & 0.5223E-02  & 0.5100E-03\\
2.3 & 0.2664E-03 & 0.2507E-04 & 0.5708E-02  & 0.6290E-04\\
4.5 & 0.8259E-02 & 0.4510E-03 & 0.4430E-01  & 0.5291E-03\\
5.0 & 0.5531E-02 & 0.7410E-03 & 0.3548E-01  & 0.8302E-03\\
\noalign{\smallskip}\hline
\end{tabular}
\end{table}
\bigskip
\bigskip
\begin{table}[htp]
\centering
\caption{CPU comparisons of Example 2.}
\label{tab:2}       
\begin{tabular}{ll}
\hline\noalign{\smallskip}
$N$ & LCM  \\
\noalign{\smallskip}\hline\noalign{\smallskip}
3 & 1.140 \\
4 & 1.258 \\
\noalign{\smallskip}\hline
\end{tabular}
\end{table}
\bigskip
\bigskip
\section{Conclusion}
\label{sec:6}
\bigskip
In this study, numerical investigation of a class of FDEs has been considered. Besides, a technique based on matrices together with collocation points has been introduced. With the help of computer programmes, Maple and Matlab, numerical approximation of Laguerre collocation method has been compared with the techniques such as Runge-Kutta and Hermite collocation methods. Moreover, the technique together with an error analysis have been implemented on some illustrations to see the applicability as well as efficiency of the method. Then the results have been seen by figures and tables. We have the advantages of the method such as straightforward applicability and coding come along with the proven convergency by some authors \cite{Khader}, \cite{Muroya}, \cite{Yang}. As a future plan, this technique can be improved and applied on some other mathematical models after completing some modifications.
\bigskip



\medskip
Received xxxx 20xx; revised xxxx 20xx.
\medskip

\end{document}